\documentstyle{article}

\newtheorem{thm}{Theorem}[section]
\newtheorem{cor}[thm]{Corollary}

\newtheorem{prop}[thm]{Proposition}
\newtheorem{defn}[thm]{Definition}
\newtheorem{rem}[thm]{Remark}

\begin{document}

\baselineskip = 15 pt
\parskip 4pt

\newtheorem{defin}{Definition}
\newcommand{\Al}{${\cal A}_{\hbar,\eta}(\widehat{sl}_2)~$}
\newcommand{\Alo}{${\cal A}_{\hbar,0}(\widehat{g})~$}
\newcommand{\Apq}{${\cal A}_{q,p}(\widehat{sl}_2)~$}
\newcommand{\Apqn}{${\cal A}_{q,p}(\widehat{sl}_N)~$}
\newcommand{\Apqg}{${\cal A}_{q,p}(\widehat{g})~$}
\newcommand{\Eg}{${\cal E}_{q,p}(\widehat{g})~$}
\newcommand{\Egc}[2]{${\cal A}^{#1}(\{\Psi\},\widehat{g})_{#2}~$}
\newcommand{\Eegc}[2]{{\cal A}^{#1}(\{\Psi\},\widehat{g})_{#2}}
\newcommand{\Alg}{${\cal A}_{\hbar,\eta}(\widehat{g})~$}
\newcommand{\Algi}[2]{{\cal A}_{\hbar,#1}(\widehat{g})_{#2}}
\newcommand{\Alp}[2]{{\cal A}_{#1,p}(\widehat{sl_2})_{#2}}
\newcommand{\Algp}[2]{{\cal A}_{#1,p}(\widehat{g})_{#2}}
\newcommand{\AlN}{${\cal A}_{\hbar,\eta}(\widehat{sl}_N)~$}
\newcommand{\phq}[4]{\phi^{(#1)}\left(\frac{#2}{#3}{#4}\right)}
\newcommand{\Phq}[4]{\Phi^{(#1)}_{ij}\left(\frac{#2}{#3}{#4}\right)}
\newcommand{\Psq}[4]{\Psi^{(#1)}\left(\frac{#2}{#3}{#4}\right)}


\centerline{\LARGE \bf Generalized quantum current algebras}

\vspace{0.5cm}

\centerline{
Liu Zhao$^{ab}$ 
}

\centerline{${}^a$ The Abdus Salam International Centre for Theoretical Physical,}
\centerline{P.O.Box 586, 34100 Trieste, Italy}
\vspace{0.2cm}
\centerline{${}^b$Institute of Modern Physics, Northwest University, Xian 710069,
China\footnote{E-mail: lzhao@nwu.edu.cn}
}




\begin{abstract}
Two general families of new quantum deformed current algebras are proposed
and identified both as infinite Hopf family of algebras, a structure which 
enable one to define ``tensor products'' of these algebras. The standard
quantum affine algebras turn out to be a very special case of both algebra
families, in which case the infinite Hopf family structure degenerates into
standard Hopf algebras. The relationship between the two algebra families 
as well as their
various special examples are discussed, and the free boson representation is
also considered. 
\end{abstract}


\section*{}
Quantum groups, since proposed by Drinfeld \cite{D1,D2}, have attracted 
the attentions of
both mathematicians and theoretical physicists for more than ten years.
And even after so many years' extensive studies, the interest in quantum groups
and their various extensions--quantum affine algebras \cite{D3,RS}, Yangian doubles
\cite{KT} and so on--is still not faded.

Theoretical physicists are fascinated about quantum groups and their extensions
because these algebraic objects are the right candidates for describing the dynamical
symmetries in certain models in integrable quantum field theories and/or exactly
solvable statistical physics. Mathematicians are interested in quantum groups
and their extensions because these are the first known nontrivial Hopf
algebras--non-commutative non-co-commutative associative algebras expected
for quite some years but not realized until the discovery of Drinfeld.

Recently, accompanying the search and investigation for elliptic quantum groups 
\cite{felder1,felder2,felder3},
it is realized that there are more general algebraic structures which
are of interests to both mathematicians and physicists. One such example is
the elliptic quantum groups proposed by Felder et al \cite{felder1,felder2,felder3}
which belong to the class
of quasi-triangular quasi-Hopf algebras \cite{D4}--certain twists \cite{jimbot}
of the standard Hopf algebra
structures. From a pure mathematical point of view, the most important
significance for the discovery of elliptic quantum groups might be that it reveals
the possibility for the existence of non-co-associative and
non-commutative non-co-commutative associative algebras. Such algebras also happen
to describe the dynamical symmetry of some statistical models and hence greatly
attracted the attention of theoretical physicists.

Further investigations along this line showed that there are several
kinds of elliptic quantum groups 
\cite{felder1,felder2,felder3,Foda1,Foda2,HY,hou2,jimbot,Konno,Z1}, 
some are recognized as quasi-triangular quasi-Hopf
algebras and some are still not. Due to the lack of co-algebraic structures
for some of the elliptic quantum groups, it is nature to
ask if the structure of quasi-triangular quasi-Hopf algebras is still adequate
for describing the co-algebraic structures for these algebras. The answer to this
question remains open and it seems very hard to get it solved.
Alternatively, it may be reasonable to try some other formulations to describe
the co-algebraic structures for some of the new elliptic quantum groups, even if
the answer to the above question finally proves to be a positive ``yes''. The infinite
Hopf family of algebras described in \cite{f, hou2} is just such a candidate.

For standard quantum affine algebras and Yangian doubles, three major realizations
are proved to exists, i.e. the Yang-Baxter type realization (also known as
Reshetikhin-Semenov-Tian-Shansky realization \cite{RS}), Drinfeld realization 
\cite{D1,D2} and Drinfeld
current realization (or simply current realization) \cite{D3}. The first and 
last realizations
are usually used to describe the co-algebraic structures and for the investigation of
infinite dimensional representations respectively. The standard Hopf algebra structure
cannot be written closely using the current realization only. However, there exists
an alternative Hopf algebra structure for quantum affine algebras and Yangian doubles
which closes over the currents only. Such a structure is first discovered by Drinfeld
in the case of quantum affine algebras and is called Drinfeld's Hopf structures \cite{Ding1}.
It is generally believed that the three realizations exist for all
quantum deformations of the universal enveloping algebras of the classical affine
Lie algebras, although only one or two realizations for some recently discovered
quantum algebras are known to exist.

In this paper, we shall study two general families of quantum current algebras which
contain many well-known quantum current algebras as special cases. We shall show that
both families of current algebras have the structure of infinite Hopf family of algebras,
which ensure a proper definition of fused representations (i.e. tensor product
representations via comultiplications).

Our investigation will be restricted in the
current realization only. Strictly speaking, the algebras we shall consider
are not algebras but ``functional algebras" \cite{Ding2, SK1} 
generalizing the concept of a usual
algebra. However we shall abuse the language in the context and call the object
under investigation ``algebras''.

Throughout the following, $g$ will be a classical simply-laced Lie algebra with
Cartan matrix $(A_{ij})$.

\section{}

\begin{defn}
Let $q,p$ be generic parameters.
We choose the set of analytic functions $\Psi_{ij}(z|q)$ ($i,j=1,...,\mbox{rank}(g)$)
such that they depend on the suffices $i,j$ only through the Cartan matrix elements
$A_{ij}$ of $g$ and that

\begin{eqnarray}
\Psi_{ij}(z|q) = \Psi_{ji}(z^{-1}|q)^{-1}.  \label{AA}
\end{eqnarray}

\noindent The functional current algebra \Egc{(1)}{q,p,c} is defined as the associative
algebra generated by the currents $E_i(z)$, $F_i(z)$, invertible $H^\pm_i(z)$ with
$i=1, 2, ..., \mbox{rank} (g)$, central element $c$ and the unit 1 with relations

\begin{eqnarray}
H^\pm_i(z)H^\pm_j(w)&=&
\frac{\Psi_{ij}\left(\left.\frac{z}{w}\right| q \right)}
{\Psi_{ij}\left(\left.\frac{z}{w}\right| \tilde{q} \right)}
H^\pm_j(w)H^\pm_i(z), \label{31}\\
H^+_i(z)H^-_j(w)&=&
\frac{\Psi_{ij}\left(\left.\frac{z}{w} p^{c/2} \right| q \right)}
{\Psi_{ij}\left(\left.\frac{z}{w} p^{-c/2} \right| \tilde{q} \right)}
H^-_j(w)H^+_i(z), \\
H^\pm_i(z)E_j(w)&=&
\Psi_{ij}\left(\left.\frac{z}{w} p^{\pm c/4} \right| q \right)
E_j(w)H^\pm_i(z),\\
H^\pm_i(z)F_j(w)&=&
\Psi_{ij}\left(\left.\frac{z}{w} p^{\mp c/4} \right| \tilde{q} \right)^{-1}
F_j(w)H^\pm_i(z),\\
E_i(z)E_j(w)&=&
\Psi_{ij}\left(\left.\frac{z}{w}\right| q \right)
E_j(w)E_i(z),\\
F_i(z)F_j(w)&=&
\Psi_{ij}\left(\left.\frac{z}{w} \right| \tilde{q} \right)^{-1}
F_j(w)F_i(z),\\
{}[ E_i(z), F_j(w) ] &=&\frac{\delta_{ij}}{(p^{\frac{1}{2}}-p^{-\frac{1}{2}})} \left[
\delta\left(\frac{z}{w}p^{-c/2}\right)H^+_i(wp^{c/4})
-\delta\left(\frac{w}{z}p^{-c/2}\right)H^-_i(zp^{c/4}) \right], \label{delta1}\\
E_i(z_1)E_i(z_2)E_j(w) \!\!\!\!&-&\!\!\!\! f^{+}_{ij}(z_1/w,z_2/w |q)
E_i(z_1)E_j(w)E_i(z_2) + E_j(w)E_i(z_1)E_i(z_2)\nonumber\\
& & + (\mbox{replacement } ~z_1 \leftrightarrow z_2) =0,
\hspace{0.5cm} A_{ij}=-1,\\
F_i(z_1)F_i(z_2)F_j(w) \!\!\!\!&-&\!\!\!\! f^{-}_{ij}(z_1/w,z_2/w| \tilde{q})
F_i(z_1)F_j(w)F_i(z_2) + F_j(w)F_i(z_1)F_i(z_2) \nonumber\\
& & + (\mbox{replacement } ~z_1 \leftrightarrow z_2) =0,
\hspace{0.5cm} A_{ij}=-1, \label{50}
\end{eqnarray}

\noindent where

\begin{eqnarray*} 
& &\tilde{q} \equiv qp^{c},\\
& &f_{ij}^+ (z_1/w,z_2/w| q) =
\frac{\left(\Psi_{ii}\left(\left.\frac{z_2}{z_1}\right| q\right)+1 \right)
\left(\Psi_{ij}\left(\left.\frac{w}{z_1}\right| q\right)
\Psi_{ij}\left(\left.\frac{w}{z_2}\right| q\right)+1 \right)}
{\Psi_{ij}\left(\left.\frac{w}{z_2}\right| q\right)
+\Psi_{ii}\left(\left.\frac{z_2}{z_1}\right| q\right)
\Psi_{ij}\left(\left.\frac{w}{z_1}\right| q\right)},\\
& &f_{ij}^- (z_1/w,z_2/w| \tilde{q}) =
\frac{\left(\Psi_{ii}\left(\left.\frac{z_2}{z_1}\right| \tilde{q}\right)^{-1}+1 \right)
\left(\Psi_{ij}\left(\left.\frac{w}{z_1}\right| \tilde{q}\right)^{-1}
\Psi_{ij}\left(\left.\frac{w}{z_2}\right| \tilde{q}\right)^{-1}+1 \right)}
{\Psi_{ij}\left(\left.\frac{w}{z_2}\right| \tilde{q}\right)^{-1}
+\Psi_{ii}\left(\left.\frac{z_2}{z_1}\right| \tilde{q}\right)^{-1}
\Psi_{ij}\left(\left.\frac{w}{z_1}\right| \tilde{q}\right)^{-1}}.
\end{eqnarray*}
\end{defn}

Let $\{ {\cal A}_n,~n\in Z\}$ be a family of associative algebras over $C$ with unit.
Let $\{v^{(n)}_i,~i=1,~...,~\mbox{dim}({\cal A}_n)\}$ be a basis of ${\cal A}_n$.
The maps

\begin{eqnarray*}
\tau_n^{\pm}: {\cal A}_n &\rightarrow& {\cal A}_{n \pm 1}\\
v^{(n)}_i &\mapsto& v^{(n\pm 1)}_i
\end{eqnarray*}

\noindent are morphisms from ${\cal A}_n$ to ${\cal A}_{n \pm 1}$.
For any two integers $n,~m$ with $n<m$, we can specify a pair of morphisms

\begin{eqnarray}
& &\tau^{(m,n)} =Mor({\cal A}_m,~{\cal A}_n) \equiv \tau_{m-1}^+...\tau_{n+1}^+\tau_n^+:~~
{\cal A}_n \rightarrow {\cal A}_m, \nonumber\\
& &\tau^{(n,m)} =Mor({\cal A}_n,~{\cal A}_m) \equiv \tau_{n+1}^-...\tau_{m-1}^-\tau_m^-:~~
{\cal A}_m \rightarrow {\cal A}_n   \label{taumn}
\end{eqnarray}

\noindent with $\tau^{(m,n)}\tau^{(n,m)}=id_{m},~\tau^{(n,m)}\tau^{(m,n)}=id_{n}$.
Clearly the morphisms $\tau^{(m,n)},~n,m\in Z$ satisfy the associativity
condition $\tau^{(m,p)}\tau^{(p,n)} =\tau^{(m,n)}$ and thus make the family of algebras
$\{ {\cal A}_n,~n\in Z\}$ into a category.

\begin{defn}
The category of algebras $\{ {\cal A}_n,~\{\tau^{(n,m)}\},~n,m\in Z\}$ is called an
infinite Hopf family of algebras if on each object ${\cal A}_n$ of the category
one can define the morphisms $\Delta^\pm_n: {\cal A}_n \rightarrow {\cal A}_n
\otimes {\cal A}_{n\pm 1}$, $\epsilon_n: {\cal A}_n \rightarrow C$ and antimorphisms
$S^\pm_n: {\cal A}_n \rightarrow {\cal A}_{n\pm 1}$ such that the following axioms hold,

\begin{itemize}
\item $(\epsilon_n \otimes id_{n+1}) \circ \Delta_n^+ = \tau_n^+,~
(id_{n-1} \otimes  \epsilon_n ) \circ \Delta_n^- = \tau_n^-$ \hfill{(a1)}

\item $m_{n+1} \circ (S^+_n \otimes id_{n+1}) \circ \Delta_n^+
= \epsilon_{n+1} \circ \tau^+_n,~
m_{n-1} \circ (id_{n-1} \otimes S^-_n) \circ \Delta_n^-
= \epsilon_{n-1} \circ \tau^-_n $ \hfill{(a2)}

\item $(\Delta_n^- \otimes id_{n+1}) \circ \Delta_n^+ =
(id_{n-1} \otimes \Delta_n^+ ) \circ \Delta_n^-$ \hfill{(a3)}
\end{itemize}

\noindent in which $m_{n}$ is the algebra multiplication for ${\cal A}_n$.
\end{defn}

Now let us consider the structure of infinite Hopf family of algebras for our
algebra \Egc{(1)}{q,p,c}. This algebra is determined uniquely by the defining relations
provided the following data are fixed: $g,~q,~p,~c$.

In general, given a series of $c_n,~n\in Z$, we can define
the series of parameters $\{ q_n, n\in Z\}$ by
the relations

\begin{eqnarray*}
q_{n+1}/q_{n}=p^{c_n},
\end{eqnarray*}

\noindent starting from the data $q_1=q$, $c_1=c$. It is obvious that
$\tilde{q}=q_{2}$. In the same fashion we can define $\tilde{q}_{n}=q_{n+1}$
and write the algebras ${\cal A}_n \equiv$ \Egc{(1)}{q_n,p,c_n}.
The generating currents $H_i^\pm(z),~E_i(z)$ and $F_i(z)$
for the algebra ${\cal A}_n$ are denoted as
$H_i^\pm(z|n),~E_i(z|n)$ and $F_i(z|n)$ respectively.

We collect the families of algebras $A=\{{\cal A}_n, n\in Z\}$.
This family of algebras can be turned into a category if we introduce
the morphisms $\tau_n^\pm$

\begin{eqnarray*}
\tau_n^\pm : {\cal A}_n &\rightarrow& {\cal A}_{n \pm 1} \\
H^\pm_i(z|n) &\mapsto& H^\pm_i(z|n\pm 1)\\
E_i(z|n) &\mapsto& E_i(z|n\pm 1)\\
F_i(z|n) &\mapsto& F_i(z|n\pm 1)\\
c_n &\mapsto& c_{n\pm 1}
\end{eqnarray*}

\noindent and defining the compositions $\tau^{(n,m)}$ as did in (\ref{taumn}).

The following proposition say that $A$ is an infinite Hopf family of algebras.

\begin{prop}
The category $A$ of algebras $\{{\cal A}_n, n\in Z\}$ form an infinite Hopf
family of algebras with the Hopf family structures given as follows:

\begin{itemize}
\item the comultiplications $\Delta_n^\pm$:
\begin{eqnarray*}
\Delta_n^+ c_n \!\!\!\!&=&\!\!\!\! c_n + c_{n+1}, \label{Co} \\
\Delta_n^+ H_i^+(z|n) \!\!\!\!&=&\!\!\!\!
H_i^+(z p^{c_{n+1}/4}|n)
\otimes H_i^+(z p^{-c_{n}/4}|n+1),\\
\Delta_n^+ H_i^-(z|n) \!\!\!\!&=&\!\!\!\!
H_i^-(z p^{-c_{n+1}/4}| n)
\otimes H_i^-(z p^{c_{n}/4}| n+1),\\
\Delta_n^+ E_i(z|n) \!\!\!\!&=&\!\!\!\! E_i(z| n)
\otimes 1 + H^-_i(z p^{c_n/4}| n) \otimes
E_i(z p^{c_n/2}| n+1),\\
\Delta_n^+ F_i(z|n) \!\!\!\!&=&\!\!\!\! 1 \otimes
F_i(z| n+1) + F_i(z p^{c_{n+1}/2}| n)
\otimes H^+_i(z p^{c_{n+1}/4}| n+1),\\
\\
\Delta_n^- c_n \!\!\!\!&=&\!\!\!\! c_{n-1} + c_n,\\
\Delta_n^- H_i^+(z|n) \!\!\!\!&=&\!\!\!\!
H_i^+(z p^{c_{n}/4}| n-1)
\otimes H_i^+(z p^{-c_{n-1}/4}| n),\\
\Delta_n^- H_i^-(z|n) \!\!\!\!&=&\!\!\!\!
H_i^-(z p^{-c_{n}/4}| n-1)
\otimes H_i^-(z p^{c_{n-1}/4}| n),\\
\Delta_n^- E_i(z|n) \!\!\!\!&=&\!\!\!\! E_i(z| n-1)
\otimes 1 + H^-_i(z p^{c_{n-1}/4}| n-1) \otimes
E_i(z p^{c_{n-1}/2}| n),\\
\Delta_n^+ F_i(z|n) \!\!\!\!&=&\!\!\!\! 1 \otimes
F_i(z| n) + F_i(z p^{c_{n}/2}| n-1)
\otimes H^+_i(z p^{c_{n}/4}| n);
\end{eqnarray*}

\item the counits $\epsilon_n$:
\begin{eqnarray*}
\epsilon_n ( c_n )\!\!\!\!&=&\!\!\!\!0,\\
\epsilon_n ( 1_n ) \!\!\!\!&=&\!\!\!\! 1,\\
\epsilon_n ( H^\pm_i(z| n))  \!\!\!\!&=&\!\!\!\! 1,\\
\epsilon_n ( E_i(z| n)) \!\!\!\! &=& \!\!\!\! 0,\\
\epsilon_n ( F_i(z| n)) \!\!\!\! &=& \!\!\!\! 0;
\end{eqnarray*}

\item the antipodes $S_n^\pm$:

\begin{eqnarray*}
S_n^\pm c_n \!\!\!\!&=&\!\!\!\! - c_{n \pm 1},\\
S_n^\pm H^+_i(z| n)\!\!\!\!& =&\!\!\!\!
[ H^+_i(z| n \pm 1 )]^{-1},\\
S_n^\pm H^-_i(z| n)\!\!\!\!& =&\!\!\!\!
[ H^-_i(z| n \pm 1 )]^{-1},\\
S^\pm_n E_i(z|n) \!\!\!\!&=&\!\!\!\!
- H^-_i(z p^{-c_{n \pm 1}/4}| n \pm 1)^{-1}
E_i(z p^{-c_{n \pm 1}/2}| n \pm 1),\\
S^\pm_n F_i(z| {n}) \!\!\!\!&=&\!\!\!\!
- F_i(z p^{-c_{n \pm 1}/2}| n \pm 1)
H^+_i(z p^{-c_{n \pm 1}/4}| n \pm 1)^{-1}.
\end{eqnarray*}
\end{itemize}
\end{prop}

Despite these unusual co-algebraic structures, we can still define tensor
product homomorphisms for both algebra families. We have

\begin{prop}
The comultiplication $\Delta^+_n$ induces an algebra homomorphism

\begin{eqnarray*}
\rho: \Eegc{(1)}{q_n,p,c_n+c_{n+1}} &\rightarrow& \Eegc{(1)}{q_n,p,c_n} \otimes
\Eegc{(1)}{q_{n+1},p,c_{n+1}}\\
X &\mapsto& \Delta^+_n \tilde{X},
\end{eqnarray*}

\noindent where $X \in \Eegc{(1)}{q_n,p,c_n+c_{n+1}},~~\tilde{X} \in
\Eegc{(1)}{q_n,p,c_n}$ and

\begin{eqnarray*}
\tilde{X}=\left\{
\begin{array}{l}
$$c_n$$\\
$$H^\pm_i(z|n)$$\\
$$E_i(z|n)$$\\
$$F_i(z|n)$$
\end{array}
\right.
\quad \quad \quad
\mbox{if}
\quad \quad \quad
X=\left\{
\begin{array}{l}
$$c_n+c_{n+1}$$\\
$$H^\pm_i(z|n)$$\\
$$E_i(z|n)$$\\
$$F_i(z|n)$$
\end{array}
\right. .
\end{eqnarray*}

\noindent Likewise, the comultiplication $\Delta^-_n$
induces an algebra homomorphism

\begin{eqnarray*}
\bar{\rho}: \Eegc{(1)}{q_{n-1},p,c_{n-1}+c_{n}} &\rightarrow&
\Eegc{(1)}{q_{n-1},p,c_{n-1}}
\otimes \Eegc{(1)}{q_n,p,c_{n}}\\
X &\mapsto& \Delta^+_n \tilde{X},
\end{eqnarray*}

\noindent where $X \in \Eegc{(1)}{q_{n-1},p,c_{n-1}+c_{n}},~~
\tilde{X} \in \Eegc{(1)}{q_n,p,c_n}$ and

\begin{eqnarray*}
\tilde{X}=\left\{
\begin{array}{l}
$$c_{n}$$\\
$$H^\pm_i(z|n)$$\\
$$E_i(z|n)$$\\
$$F_i(z|n)$$
\end{array}
\right.
\quad \quad \quad
\mbox{if}
\quad \quad \quad
X=\left\{
\begin{array}{l}
$$c_{n-1}+c_{n}$$\\
$$H^\pm_i(z|n-1)$$\\
$$E_i(z|n-1)$$\\
$$F_i(z|n-1)$$
\end{array}
\right. .
\end{eqnarray*}
\end{prop}

\begin{cor}
Let $m$ be a positive integer. The iterated comultiplication
$\Delta^{(m)+}_n=(id_n \otimes id_{n+1} \otimes ...
\otimes id_{n+m-2} \otimes \Delta_{n+m-1}^+)
\circ (id_n \otimes id_{n+1} \otimes ...
\otimes id_{n+m-3} \otimes \Delta_{n+m-2}^+) ...
\circ (id_n \otimes \Delta_{n+1}^+) \circ \Delta_n^+$ induces an algebra homomorphism
$\rho^{(m)}$

\begin{eqnarray*}
\rho^{(m)}: & &\Eegc{(1)}{q^{(n)},p,c_n +c_{n+1} + ... + c_{n+m}}\\
& & \rightarrow
\Eegc{(1)}{q^{(n)},p,c_n} \otimes \Eegc{1}{q^{(n+1)},c_{n+1}} \otimes ... \otimes
\Eegc{(1)}{q^{(n+m)},p,c_{n+m}}
\end{eqnarray*}

\noindent in the spirit of Proposition 1.4.
\end{cor}

\begin{rem}
We stress here that the maps $\rho,~\bar{\rho}$ and $\rho^{(m)}$ are {\em algebra
homomorphisms}, whilst $\tau_n^\pm,~\tau^{(n,m)}$ and $\Delta_n^\pm$ etc are
only algebra morphisms. The difference between algebra morphisms and algebra
homomorphisms lies in that the latter preserves the structure functions whilst the
former does not.
\end{rem}

\section{}

\begin{defn}
Let $q,\tilde{q}$ be generic parameters.
We choose the set of analytic functions $\Psi_{ij}(z|q)$ ($i,j=1,...,\mbox{rank}(g)$)
such that they obey the condition (\ref{AA}).
The functional current algebra \Egc{(2)}{q,\tilde{q},\beta,\gamma}
is defined as the associative
algebra generated by the currents $E_i(z)$, $F_i(z)$, invertible $H^\pm_i(z)$ with
$i=1, 2, ..., \mbox{rank} (g)$, central elements $\beta, \gamma$,
and the unit 1 with relations

\begin{eqnarray}
H^\pm_i(z)H^\pm_j(w)&=&
\frac{\Psi_{ij}\left(\left.\frac{z}{w}\right| q \right)}
{\Psi_{ij}\left(\left.\frac{z}{w}\right| \tilde{q} \right)}
H^\pm_j(w)H^\pm_i(z), \label{36}\\
H^+_i(z)H^-_j(w)&=&
\frac{\Psi_{ij}\left(\left.\frac{z}{w} (\beta^{-1}\gamma)^{1/2} \right| q \right)}
{\Psi_{ij}\left(\left.\frac{z}{w} (\beta^{-1}\gamma)^{-1/2} \right| \tilde{q} \right)}
H^-_j(w)H^+_i(z), \\
H^+_i(z)E_j(w)&=&
\Psi_{ij}\left(\left.\frac{z}{w} \beta^{-1/2} \right| q \right)
E_j(w)H^+_i(z),\\
H^-_i(z)E_j(w)&=&
\Psi_{ij}\left(\left.\frac{z}{w} \gamma^{-1/2} \right| q \right)
E_j(w)H^-_i(z),\\
H^+_i(z)F_j(w)&=&
\Psi_{ij}\left(\left.\frac{z}{w} \beta^{1/2} \right| \tilde{q} \right)^{-1}
F_j(w)H^+_i(z),\\
H^-_i(z)F_j(w)&=&
\Psi_{ij}\left(\left.\frac{z}{w} \gamma^{1/2} \right| \tilde{q} \right)^{-1}
F_j(w)H^-_i(z),\\
E_i(z)E_j(w)&=&
\Psi_{ij}\left(\left.\frac{z}{w}\right| q \right)
E_j(w)E_i(z),\\
F_i(z)F_j(w)&=&
\Psi_{ij}\left(\left.\frac{z}{w} \right| \tilde{q} \right)^{-1}
F_j(w)F_i(z),\\
{}[ E_i(z), F_j(w) ] &=&\frac{\delta_{ij}}{(\tilde{q}/q-1)} \left[
\delta\left(\frac{z}{w}\beta \right)H^+_i(w \beta^{-1/2})
-\delta\left(\frac{w}{z} \gamma^{-1} \right) H^-_i(z \gamma^{1/2}) \right], \label{delta2}\\
E_i(z_1)E_i(z_2)E_j(w) \!\!\!\!&-&\!\!\!\! f^{+}_{ij}(z_1/w,z_2/w |q)
E_i(z_1)E_j(w)E_i(z_2) + E_j(w)E_i(z_1)E_i(z_2)\nonumber\\
& & + (\mbox{replacement } ~z_1 \leftrightarrow z_2) =0,
\hspace{0.5cm} A_{ij}=-1,\\
F_i(z_1)F_i(z_2)F_j(w) \!\!\!\!&-&\!\!\!\! f^{-}_{ij}(z_1/w,z_2/w| \tilde{q})
F_i(z_1)F_j(w)F_i(z_2) + F_j(w)F_i(z_1)F_i(z_2) \nonumber\\
& & + (\mbox{replacement } ~z_1 \leftrightarrow z_2) =0,
\hspace{0.5cm} A_{ij}=-1, \label{61}
\end{eqnarray}

\noindent where

\begin{eqnarray*} 
& &f_{ij}^+ (z_1/w,z_2/w| q) =
\frac{\left(\Psi_{ii}\left(\left.\frac{z_2}{z_1}\right| q\right)+1 \right)
\left(\Psi_{ij}\left(\left.\frac{w}{z_1}\right| q\right)
\Psi_{ij}\left(\left.\frac{w}{z_2}\right| q\right)+1 \right)}
{\Psi_{ij}\left(\left.\frac{w}{z_2}\right| q\right)
+\Psi_{ii}\left(\left.\frac{z_2}{z_1}\right| q\right)
\Psi_{ij}\left(\left.\frac{w}{z_1}\right| q\right)},\\
& &f_{ij}^- (z_1/w,z_2/w| \tilde{q}) =
\frac{\left(\Psi_{ii}\left(\left.\frac{z_2}{z_1}\right| \tilde{q}\right)^{-1}+1 \right)
\left(\Psi_{ij}\left(\left.\frac{w}{z_1}\right| \tilde{q}\right)^{-1}
\Psi_{ij}\left(\left.\frac{w}{z_2}\right| \tilde{q}\right)^{-1}+1 \right)}
{\Psi_{ij}\left(\left.\frac{w}{z_2}\right| \tilde{q}\right)^{-1}
+\Psi_{ii}\left(\left.\frac{z_2}{z_1}\right| \tilde{q}\right)^{-1}
\Psi_{ij}\left(\left.\frac{w}{z_1}\right| \tilde{q}\right)^{-1}}.
\end{eqnarray*}
\end{defn}

This algebra is determined uniquely by the defining relations
provided the following data are fixed: $g,~q,~\tilde{q}$.

Now let us choose an arbitrary set of parameters
$\{q^{(n)}, ~n\in Z\}$ and define $\tilde{q}^{(n)}=q^{(n+1)}$.
We collect the family of algebras
$B=\{ {\cal B}_n, ~n\in Z\}$ where
${\cal B}_n \equiv \Eegc{(2)}{q^{(n)},\tilde{q}^{(n)},\beta_n,\gamma_n}$.
The generating currents $H_i^\pm(z),~E_i(z)$ and $F_i(z)$
for the algebra ${\cal B}_n$ are denoted as
$H_i^\pm(z|n),~E_i(z|n)$ and $F_i(z|n)$ respectively.

The algebra family $B$ can also be turned into a category in the same way as
we did for the family $A$, provided the basic morphisms $\tau_n^\pm$ are given
as follows,

\begin{eqnarray*}
\tau_n^\pm : {\cal B}_n &\rightarrow& {\cal B}_{n \pm 1} \\
\beta_n &\mapsto& \beta_{n\pm 1}\\
\gamma_n &\mapsto& \gamma_{n\pm 1}\\
H^\pm_i(z|n) &\mapsto& H^\pm_i(z|n\pm 1)\\
E_i(z|n) &\mapsto& E_i(z|n\pm 1)\\
F_i(z|n) &\mapsto& F_i(z|n\pm 1).
\end{eqnarray*}

The following proposition say that the algebra family $B$
is also an infinite Hopf family of algebras.

\begin{prop}
The category $B$ of algebras $\{{\cal B}_n, n\in Z\}$ form an infinite Hopf
family of algebras with the Hopf family structures given as follows:

\begin{itemize}
\item the comultiplications $\Delta_n^\pm$:
\begin{eqnarray*}
\Delta_n^+ \beta_n \!\!\!\!&=&\!\!\!\! \beta_n\beta_{n+1}, \label{CB} \\
\Delta_n^+ \gamma_{n} \!\!\!\!&=&\!\!\!\!
\gamma_{n} \gamma_{n+1},\\
\Delta_n^+ H_i^+(z|n) \!\!\!\!&=&\!\!\!\!
H_i^+(z (\beta_{n+1})^{-1/2}| n)
\otimes H_i^+(z (\beta_{n})^{1/2}| n+1),\\
\Delta_n^+ H_i^-(z|n) \!\!\!\!&=&\!\!\!\!
H_i^-(z (\gamma_{n+1})^{-1/2}| n)
\otimes H_i^-(z (\gamma_n)^{1/2}| n+1),\\
\Delta_n^+ E_i(z|n) \!\!\!\!&=&\!\!\!\! E_i(z|n)
\otimes 1 + H^-_i(z (\gamma_n)^{1/2}| n) \otimes
E_i(z \gamma_n| n+1),\\
\Delta_n^+ F_i(z|n) \!\!\!\!&=&\!\!\!\!
1 \otimes F_i(z| n+1) + F_i(z (\beta_{n+1})^{-1}|n)
\otimes H^+_i(z (\beta_{n+1})^{-1/2}| n+1),\\
\\
\Delta_n^- \beta_n \!\!\!\!&=&\!\!\!\! \beta_{n-1}\beta_{n},  \\
\Delta_n^- \gamma_{n} \!\!\!\!&=&\!\!\!\!
\gamma_{n-1} \gamma_{n},\\
\Delta_n^- H_i^+(z|n) \!\!\!\!&=&\!\!\!\!
H_i^+(z (\beta_n)^{-1/2}|n-1) \otimes H_i^+(z (\beta_{n-1})^{1/2}|n),\\
\Delta_n^- H_i^-(z|n) \!\!\!\!&=&\!\!\!\!
H_i^-(z (\gamma_n)^{-1/2}|n-1) \otimes H_i^-(z (\gamma_{n-1})^{1/2}|n),\\
\Delta_n^- E_i(z|n) \!\!\!\!&=&\!\!\!\!
E_i(z|n-1)\otimes 1 +H_i^{-}(z(\gamma_{n-1})^{1/2}|n-1) \otimes
E_i(z\gamma_{n-1}|n),\\
\Delta_n^- F_i(z|n) \!\!\!\!&=&\!\!\!\!
1 \otimes F_i(z| n) + F_i(z (\beta_{n})^{-1}|n-1)
\otimes H^+_i(z (\beta_{n})^{-1/2}| n);\\
\end{eqnarray*}

\item the counits $\epsilon_n$:
\begin{eqnarray*}
\epsilon_n ( \beta_n )\!\!\!\!&=&\!\!\!\!1,\\
\epsilon_n ( \gamma_n )\!\!\!\!&=&\!\!\!\!1,\\
\epsilon_n ( 1_n ) \!\!\!\!&=&\!\!\!\! 1,\\
\epsilon_n ( H^\pm_i(z|n))  \!\!\!\!&=&\!\!\!\! 1,\\
\epsilon_n ( E_i(z|n)) \!\!\!\! &=& \!\!\!\! 0,\\
\epsilon_n ( F_i(z|n)) \!\!\!\! &=& \!\!\!\! 0;
\end{eqnarray*}

\item the antipodes $S_n^\pm$:

\begin{eqnarray*}
S_n^\pm \beta_n \!\!\!\!&=&\!\!\!\! (\beta_{n \pm 1})^{-1},\\
S_n^\pm \gamma_n \!\!\!\!&=&\!\!\!\! (\gamma_{n \pm 1})^{-1},\\
S_n^\pm H^+_i(z| n)\!\!\!\!& =&\!\!\!\!
[ H^+_i(z| n \pm 1)]^{-1},\\
S_n^\pm H^-_i(z| n)\!\!\!\!& =&\!\!\!\!
[ H^-_i(z| n \pm 1)]^{-1},\\
S^\pm_n E_i(z|n)\!\!\!\!&=&\!\!\!\!-H_i^{-}(z(\gamma_{n\pm 1})^{-1/2}|n\pm 1)^{-1}
E_i(z(\gamma_{n\pm 1})^{-1}|n\pm 1),\\
S^\pm_n F_i(z|n)\!\!\!\!&=&\!\!\!\!
-F_i(z\beta_{n\pm 1}|n\pm 1)
H_i^{+}(z(\beta_{n\pm 1})^{1/2}|n\pm 1)^{-1}.
\end{eqnarray*}
\end{itemize}
\end{prop}

The co-algebraic structure for the family $B$ is also rather unusual, and
we can understand such a structure more deeply by considering the tensor
product homomorphisms for such algebra families. We have

\begin{prop}
The comultiplication $\Delta^+_n$ induces an algebra homomorphism

\begin{eqnarray*}
& &\rho: \Eegc{(2)}{q^{(n)},q^{(n+2)},\beta_n\beta_{n+1},\gamma_n\gamma_{n+1}}
\rightarrow\\
& &~~~~~~~~~~ \Eegc{(2)}{q^{(n)},q^{(n+1)},\beta_n,\gamma_n} \otimes
\Eegc{(2)}{q^{(n+1)},q^{(n+2)},\beta_{n+1},\gamma_{n+1}}\\
& &~~~~~~~~~~ X \mapsto \Delta^+_n \tilde{X},
\end{eqnarray*}

\noindent where $X \in \Eegc{(2)}{q^{(n)},q^{(n+2)},\beta_n\beta_{n+1},
\gamma_n\gamma_{n+1}},~~\tilde{X} \in
\Eegc{(2)}{q^{(n)},q^{(n+1)},\beta_n,\gamma_n}$ and

\begin{eqnarray*}
\tilde{X}=\left\{
\begin{array}{l}
$$\beta_n$$\\
$$\gamma_n$$\\
$$H^\pm_i(z|n)$$\\
$$E_i(z|n)$$\\
$$F_i(z|n)$$
\end{array}
\right.
\quad \quad \quad
\mbox{if}
\quad \quad \quad
X=\left\{
\begin{array}{l}
$$\beta_n \beta_{n+1}$$\\
$$\gamma_n \gamma_{n+1}$$\\
$$H^\pm_i(z|n)$$\\
$$E_i(z|n)$$\\
$$F_i(z|n)$$
\end{array}
\right. .
\end{eqnarray*}

\noindent Likewise, the comultiplication $\Delta^-_n$
induces an algebra homomorphism

\begin{eqnarray*}
& &\bar{\rho}: \Eegc{(2)}{q^{(n-1)},q^{(n+1)},\beta_{n-1}\beta_{n},\gamma_{n-1}\gamma_{n}}
\rightarrow \\
& &~~~~~~~~~~\Eegc{(2)}{q^{(n-1)},q^{(n)},\beta_{n-1},\gamma_{n-1}}\otimes
\Eegc{(2)}{q^{(n)},q^{(n+1)},\beta_n,\gamma_n}\\
& &~~~~~~~~~~X \mapsto \Delta^+_n \tilde{X},
\end{eqnarray*}

\noindent where $X \in \Eegc{(2)}{q^{(n-1)},q^{(n+1)},\beta_{n-1}\beta_{n},
\gamma_{n-1}\gamma_{n}},~~\tilde{X} \in \Eegc{(2)}{q^{(n)},q^{(n+1)},\beta_n,
\gamma_n}$ and

\begin{eqnarray*}
\tilde{X}=\left\{
\begin{array}{l}
$$\beta_n$$\\
$$\gamma_n$$\\
$$H^\pm_i(z|n)$$\\
$$E_i(z|n)$$\\
$$F_i(z|n)$$
\end{array}
\right.
\quad \quad \quad
\mbox{if}
\quad \quad \quad
X=\left\{
\begin{array}{l}
$$\beta_{n-1}\beta_n$$\\
$$\gamma_{n-1}\gamma_n$$\\
$$H^\pm_i(z|n-1)$$\\
$$E_i(z|n-1)$$\\
$$F_i(z|n-1)$$
\end{array}
\right. .
\end{eqnarray*}
\end{prop}

\begin{cor}
Let $m$ be a positive integer. The iterated comultiplication
$\Delta^{(m)+}_n=(id_n \otimes id_{n+1} \otimes ...
\otimes id_{n+m-2} \otimes \Delta_{n+m-1}^+)
\circ (id_n \otimes id_{n+1} \otimes ...
\otimes id_{n+m-3} \otimes \Delta_{n+m-2}^+) ...
\circ (id_n \otimes \Delta_{n+1}^+) \circ \Delta_n^+$ induces an algebra homomorphism
$\rho^{(m)}$

\begin{eqnarray*}
\rho^{(m)}: & &\Eegc{(2)}{q^{(n)},q^{(n+m+1)}, \beta_n \beta_{n+1} ... \beta_{n+m},
\gamma_n \gamma_{n+1} ... \gamma_{n+m}}\\
& & \rightarrow
\Eegc{(2)}{q^{(n)},q^{(n+1)},\beta_n,\gamma_n} \otimes
\Eegc{(2)}{q^{(n+1)},q^{(n+2)},\beta_{n+1},\gamma_{n+1}} \otimes ... \\
& &~~~~\otimes
\Eegc{(2)}{q^{(n+m)},q^{(n+m+1)},\beta_{n+m},\gamma_{n+m}}
\end{eqnarray*}

\noindent in the spirit of Proposition 2.3.
\end{cor}

\begin{rem}
From the point of view of quantized affine algebras, the family $B$ of algebras
seems more natural and symmetric: with the quantized (or deformed) Cartan part of
currents splitted into positive and negative halfs, why should the central element
remain as a whole?
\end{rem}

\section{}
Now it is the point to consider the relationship between the algebras in the
families $A$ and $B$. We note that the algebra \Egc{(2)}{q,\tilde{q},\beta,\gamma}
in the family $B$ has one more generator than the algebra \Egc{(1)}{q,p,c}
in the family $A$ and hence both algebras cannot be identical in general.
However, it is possible to introduce certain restrictions to the algebras
in the family $B$ so that the algebras in the family $B$ can be
related to the one in the family $A$. That means, the family $A$ is some special
case of the family $B$.

Now we illustrate some examples of such restrictions in due course.

First let us recall that, in the algebra \Egc{(2)}{q,\tilde{q},\beta,\gamma},
the parameters $q,\tilde{q}$ are generic and their ratio is not assumed to
be related to the central elements $\beta,\gamma$. This is in
contrast to the case of \Egc{(1)}{q,p,c} in which the ratio of $q$ and $\tilde{q}$
is related to $p^c$. Now if we can consider the algebra
\Egc{(2)}{q,\tilde{q},\beta,\gamma} with the restrictions $\beta=\gamma^{-1}
=p^{-c/2}$, $\tilde{q}=qp^c$ where $p$ is some constant, then the algebra
\Egc{(2)}{q,\tilde{q},\beta,\gamma} will become \Egc{(1)}{q,p,c}, with corresponding
generating currents identified. In this special case, the co-algebraic structures of
both algebras also coincide.

Another special case is given as follows. In the algebra
\Egc{(2)}{q,\tilde{q},\beta,\gamma}, let $\tilde{q}/q=p^c$, $p$
being some constant. Then the map

\begin{eqnarray*}
\mu: \Eegc{(1)}{q,p,c} &\rightarrow& \Eegc{(2)}{q,\tilde{q},\beta,\gamma}\\
E_i(z) &\mapsto& E_i(z\gamma^{-1/2})\\
F_i(z) &\mapsto& F_i(z\beta^{1/2})\\
H^\pm_i(z) &\mapsto& H^\pm_i(z(\beta^{-1}\gamma)^{-1/4})\\
p^c &\mapsto& \beta^{-1}\gamma
\end{eqnarray*}

\noindent gives a homomorphism from \Egc{(1)}{q,p,c}
to \Egc{(2)}{q,\tilde{q},\beta,\gamma} as associative algebras.
However, under this case, the co-algebraic
structures for the algebra \Egc{(1)}{q,p,c} are not mapped into those for
the algebra \Egc{(2)}{q,\tilde{q},\beta,\gamma}.

It is interesting to mention that, the homomorphism $\mu$, though cannot map the
co-algebraic structures correctly, can provide a way of obtaining certain bosonic
realizations of the algebra \Egc{(2)}{q,\tilde{q},\beta,\gamma} from that of the
algebra \Egc{(1)}{q,p,c}, and vice versa. The simplest starting point 
will be the $c=1$ bosonic
realization for the algebra \Egc{(1)}{q,p,c}, which would lead to a bosonic realization
for the algebra \Egc{(2)}{q,\tilde{q},\beta,\gamma} at $\beta=q, \gamma=\tilde{q}=qp$.
We shall come back to this point in Section 5. Before going to the bosonic
representations, we would like to present some concrete examples for the structure
functions $\Psi_{ij}(z|q)$ to show how general our algebras are.

\section{}

In this section we present some examples for the structure functions
$\Psi_{ij}(z|q)$ of the algebras \Egc{(1)}{q,p,c} and
\Egc{(2)}{q,\tilde{q},\beta,\gamma}.

First comes a special case in which $\Psi_{ij}(z|q)$
are generic analytic functions of $z$ satisfying the condition (\ref{AA})
but are independent of $q$. Then the parameters $q_n$ will not appear at all
in both the family $A$ and $B$, and the whole family $A$
will degenerate into a single standard Hopf algebra which is nothing but the generalized
quantum current algebra given by Ding and Iohara in \cite{Ding2} (we should change the
notation $p\rightarrow q^{2}$ to compare with \cite{Ding2}). The family $B$ will also
degenerate into a single Hopf algebra which is, to our knowledge, not considered
elsewhere earlier.

Let $\psi(z)$ be an analytical function of $z$ such that

\begin{eqnarray*}
\psi(z)=-z\psi(z^{-1}),
\end{eqnarray*}

\noindent whose definition may depends on the parameters $(q,p)$ or $(q,\tilde{q})$.
We define

\begin{eqnarray}
\Psi_{ij}(z|q) = (-1)^{A_{ij}} x^{-A_{ij}}
\frac{\psi(zx^{A_{ij}})}{\psi(zx^{-A_{ij}})}, \label{BB}
\end{eqnarray}

\noindent where $x$ is an arbitrary function of the deformation parameters,
$x=x(q,p)$ for the case of \Egc{(1)}{q,p,c} and
$x=x(q,\tilde{q})$ for the case of \Egc{(2)}{q,\tilde{q},\beta,\gamma}.
The functions $\Psi_{ij}(z|q)$ given in (\ref{BB}) fulfill the condition (\ref{AA})
and hence can be used to give examples for either the algebra \Egc{(1)}{q,p,c} or
the algebra \Egc{(2)}{q,\tilde{q},\beta,\gamma}.

We note that the signature factor $(-1)^{A_{ij}}$ does not affect the
condition (\ref{AA}) and thus can be omitted from eq. (\ref{BB}), and that would of cause
lead to a slightly different definition of the algebras.

\begin{eqnarray*}
\Psi_{ij}(z|q) = x^{-A_{ij}}
\frac{\psi(zx^{A_{ij}})}{\psi(zx^{-A_{ij}})}.
\end{eqnarray*}

Let us present some more special cases for the function $\psi(z)$.

The rational function

\begin{eqnarray*}
\psi(z)=1-z
\end{eqnarray*}

\noindent is our first choice. If we further choose $x=p^{1/2}$ in this case,
the corresponding family $A$ will be identical to $U_q(\hat{g})$ with $q=p^{1/2}$.

Let $\theta_q(z)$ be the elliptic function

\begin{eqnarray*}
& & \theta_q(z) =(z|q)_\infty (qz^{-1}|q)_\infty (q| q)_\infty,\\
& &(z| q_1, ..., q_m)_\infty = \prod_{i_1,i_2,...,i_m=0}^\infty
(1-zq_1^{i_1} q_2^{i_2} ... q_m^{i_m}).
\end{eqnarray*}

\noindent We can choose

\begin{eqnarray*}
\psi(z)=\theta_y(z)
\end{eqnarray*}

\noindent with $y=y(q,p)$ for the case of \Egc{(1)}{q,p,c}
and $y=y(q,\tilde{q})$ for the case of \Egc{(2)}{q,\tilde{q},\beta,\gamma}.

If $y=const.$ and $x=p^{1/2}$, the algebra \Egc{(1)}{q,p,c} will become the
elliptic quantum group mentioned in \cite{felder1,felder2,felder3}.

It is interesting to mention that the algebra \Egc{(1)}{q,p,c} with $\psi(z)$ chosen as
$\theta_q(z)$ and $x$ chosen as $p^{1/2}$ form an elliptic generalization
of the algebra \Alg, which is the representative of the first known infinite
Hopf family of algebras. In the scaling limit, the elliptic algebra family $A$ will
become the algebra family containing \Alg. 
With the same choices of $\psi(z)$ and
$x$, the algebra \Egc{(2)}{q,\tilde{q},\beta,\gamma} will become a
new type of elliptic quantum group, which, in the case of $\beta=q, \gamma=\tilde{q}$,
will tend to the algebra of modifies screening currents with identification of parameters
$p=\tilde{q}/q$.

We emphasis that our construction enable us to introduce more free parameters and
obtain multi-parameter quantum current algebras from a very general setting.

Last we should mention that this manuscript considers only those algebras with
multiplicative spectral parameters. We could as well consider the cases with
additive spectral parameters, and in those cases the structure functions
(denoted as $\Psi_{ij}(u|\eta)$) should behave as

\begin{eqnarray*}
\Psi_{ij}(-u|\eta)=\Psi_{ji}(u|\eta)^{-1}.
\end{eqnarray*}

\noindent The algebra \Alg is actually a concrete example for the algebra family
$A$ with structure functions given by

\begin{eqnarray*}
\Psi_{ij}(u|\eta)=\frac{\mbox{sh}~\pi \eta(u-i\hbar A_{ij}/2)}
{\mbox{sh}~\pi \eta(u+i\hbar A_{ij}/2)}.
\end{eqnarray*}

\section{}

Having established the infinite Hopf family of algebras structure of
the algebra families $A$ and $B$, we now turn to consider their simplest
infinite dimensional representation, i.e. the free boson realization.

The purpose of this section is to conduct a general method for obtaining
a particular (lowest level) free boson realization for a given 
quantum current algebra. We shall show
that the free boson realization for the algebra \Egc{(1)}{q,p,c} at $c=1$
and that for the algebra \Egc{(2)}{q,\tilde{q},\beta,\gamma} at $\beta=q$,
$\gamma=\tilde{q}\equiv qp$ can actually be obtained from the same set of 
bosonic fields.

Let us start our construction by considering the generating relations for
$E_i(z) E_j(w)$, $F_i(z) F_j(w)$ and the commutator relations $[ E_i(z), F_j(w)]$
respectively. We notice that the $E_i(z) E_j(w)$, $F_i(z) F_j(w)$ relations
are the same for both algebras under consideration. This is a very important
feature for our consideration. In order to obtain a free boson realization,
we need to introduce some Riemann decomposition for the structure functions
$\Psi_{ij}(z|q)$. Suppose this decomposition is given by
\footnote{Actually, the standard definition of the $q$-affine current algebra
was given with the Riemann decomposition explicitly introduced into the
$E_i(z) E_j(w)$, $F_i(z) F_j(w)$ relations like in (\ref{Riemann1}) and
(\ref{Riemann2}).}

\begin{equation}
\Psi_{ij}(z|q)=\frac{\Phi_{ij}(z|q)}{\Phi_{ji}(z^{-1}|q)},
\end{equation}

\noindent then we can rewrite the $E_i(z) E_j(w)$, $F_i(z) F_j(w)$ relations in
the form

\begin{eqnarray}
& & \Phi_{ji}\left(\left.\frac{w}{z}\right|q\right) E_i(z) E_j(w)=
\Phi_{ij}\left(\left.\frac{z}{w}\right|q\right) E_j(w) E_i(z), \label{Riemann1}\\
& & \Phi_{ij}\left(\left.\frac{z}{w}\right|\tilde{q}\right) F_i(z) F_j(w)=
\Phi_{ji}\left(\left.\frac{w}{z}\right|\tilde{q}\right) F_j(w) F_i(z). \label{Riemann2}
\end{eqnarray}

In order to obtain a bosonic realization for the above relations, it is enough
to write down some bosonic expressions also denoted $E_i(z)$ and $F_i(z)$ such that
they satisfy the relations

\begin{eqnarray}
& &\Phi_{ji}\left(\left.\frac{w}{z}\right|q\right) E_i(z) E_j(w)
= :E_i(z) E_j(w):, \label{normal1}\\
& &\Phi_{ij}\left(\left.\frac{z}{w}\right|\tilde{q}\right) F_i(z) F_j(w)
= :F_i(z) F_j(w):,
\end{eqnarray}

\noindent where :\mbox{  }: means the standard normal ordering of bosonic expressions.
In the meantime, we set

\begin{eqnarray}
& &E_i(z) F_j(w) = \Upsilon_{ij}\left(z, w\right)
:E_i(z) F_j(w):,\\
& &F_j(w) E_i(z) = \bar{\Upsilon}_{ji}\left(z, w\right)
:E_i(z) F_j(w):, \label{normal2}
\end{eqnarray}

\noindent where, in order to yield the $\delta$-function terms in the
relations (\ref{delta1}) and (\ref{delta2}), we need the functions
$\Upsilon_{ij}\left(z, w \right)$ and $\bar{\Upsilon}_{ji}\left(z, w\right)$
to satisfy the relation

\begin{eqnarray}
\Upsilon_{ij}\left(z, w\right) -\bar{\Upsilon}_{ji}\left(z, w\right)
=\frac{\delta_{ij}}{(p^{\frac{1}{2}}-p^{-\frac{1}{2}})} \left[
\delta\left(\frac{z}{w}p^{-1/2}\right) g^{(1)+}(wp^{1/4})
-\delta\left(\frac{w}{z}p^{-1/2}\right)g^{(1)-}(zp^{1/4}) \right]   \label{CC}
\end{eqnarray}

\noindent for the case of the algebra \Egc{(1)}{q,p,c} at $c=1$ or

\begin{eqnarray*}
\Upsilon_{ij}\left(z, w\right) -\bar{\Upsilon}_{ji}\left(z, w\right)
=\frac{\delta_{ij}}{(\tilde{q}/q-1)} \left[
\delta\left(\frac{z}{w}q\right) g^{(2)+}(wq^{-1/2})
-\delta\left(\frac{w}{z} \tilde{q}^{-1}\right) g^{(2)-}(z\tilde{q}^{1/2})\right]
\end{eqnarray*}

\noindent for the case of the algebra \Egc{(2)}{q,\tilde{q},\beta,\gamma} at
$\beta=q, \gamma=\tilde{q}$, where $g^{(l)\pm}(z), l=1,2$ are some power functions of the 
arguments which could be absorbed into the definition of $H_i^\pm(z)$ as normalization factors.
Since the above two algebras at the given values of
central elements are homomorphic as mentioned in the end of Section 3, 
we shall proceed with only the
case of \Egc{(1)}{q,p,c} at $c=1$ and obtain the case of
\Egc{(2)}{q,\tilde{q},\beta,\gamma} at $\beta=q, \gamma=\tilde{q}$ as a trivial result
of the homomorphism.

From experiences in the study of bosonic realizations for standard $q$-affine algebras,
we know that the choices

\begin{eqnarray*}
& &\Upsilon_{ij}\left(z ,w\right)=\left\{
\begin{array}{ll}
    $$\frac{1}{z^2\left(1-\frac{w}{z}p^{1/2}\right)
    \left(1-\frac{w}{z}p^{-1/2}\right)}$$ & $$\mbox{for } A_{ij}=2,$$ \\
    $$\sim \mbox{some regular expressions}& $$\mbox{for } A_{ij}=-1, 0$$
\end{array}
\right.\\
& &\bar{\Upsilon}_{ji}(z, w)=\Upsilon_{ij}(w, z)
\end{eqnarray*}

\noindent fulfills the condition (\ref{CC}), with $g^{(1)+}(z)=g^{(1)-}(z)=z^{-2}$.

According to the above analysis, we now introduce the ansatz for the bosonic expressions
$E_i(z)$ and $F_i(z)$

\begin{eqnarray}
& &E_i(z) = 
:\exp \varphi_i(z): ,\\
& &F_i(z) = 
:\exp \psi_i(z): ,
\end{eqnarray}

\noindent where 

\begin{eqnarray}
& &\varphi_i(z)= Q_i + \mbox{log}(Az) P_i + \sum_{n \neq 0}  u[n] a_i[n] z^{-n} ,\\
& &\psi_i(z) = -Q_i - \mbox{log}(Bz) P_i - \sum_{n \neq 0} v[n] a_i[n] z^{-n} ,
\end{eqnarray}

\noindent $A$ and $B$ are some constants
to be related to the deformation parameters, $u[n]$ and $v[n]$ are all functions of
the integer $n$ which are independent of $z$, and $Q_i, P_i$ and $a_i[n]$ are bosonic operators
whose commutation relations are to be determined. The operators $Q_i$ and $P_i$ here play
the role of zero mode generators for the bosonic fields $\varphi_i(z)$ and $\psi_i(z)$.

Following from the above ansatz, we can write immediately

\begin{eqnarray}
& &E_i(z)E_j(w) = \mbox{exp} \langle \varphi_i(z) \varphi_j(w) \rangle
:E_i(z) E_j(w):,\\
& &F_i(z)F_j(w) = \mbox{exp} \langle \psi_i(z) \psi_j(w) \rangle
:F_i(z) F_j(w):,\\
& &E_i(z) F_j(w) = \mbox{exp} \langle \varphi_i(z) \psi_j(w) \rangle
:E_i(z) F_j(w):,
\end{eqnarray}

\noindent which, compared to eqs. (\ref{normal1}-\ref{normal2}), yield

\begin{eqnarray}
& &\langle \varphi_i(z) \varphi_j(w) \rangle 
= -\mbox{log}\left[ \Phi_{ji}\left(\left.\frac{w}{z}\right|q \right) \right],\\
& &\langle \psi_i(z) \psi_j(w) \rangle
= -\mbox{log} \left[ \Phi_{ij}\left(\left. \frac{z}{w}\right|\tilde{q}\right)\right],\\
& &\langle \varphi_i(z) \psi_j(w) \rangle = \mbox{log}\left[ \Upsilon_{ij}(z,w) \right].
\end{eqnarray}

\noindent For any concrete set of functions $\Phi_{ij}(z|q)$, these last equations serve
as a good starting point to determine the unknown coefficients $A, B, u[n], v[n]$ as well
as the unknown commutation relations for $Q_i, P_i$ and $a_i[n]$ respectively. Of cause the solution
need not to be unique.

In the concrete case when $\Psi_{ij}(z|q)$
is given by eq. (\ref{BB}) with $\psi(z)=\theta_q(z)$ and $x=p^{1/2}$, i.e.

\begin{eqnarray}
\Psi_{ij}(z|q)=(-1)^{A_{ij}}p^{-A_{ij}/2}\frac{\theta_q(zp^{A_{ij}/2})}
{\theta_q(zp^{-A_{ij}/2})}, \label{dd}
\end{eqnarray}

\noindent we have the following explicit
result. 

First the Heisenberg algebra ${\cal H}_{q,p}(g)$ with generators
$a_i[n],~P_i,~Q_i, i=1,...,$ $\mbox{rank}(g), n\in Z\backslash \{0\}$ can be introduced by
writing down the generating relations

\begin{eqnarray*}
& &[ a_i[n],~a_j[m] ] = \frac{1}{n}
\frac{(1-q^{-n})(p^{nA_{ij}/2}-p^{-nA_{ij}/2})(1-(pq)^n)}{1-p^n} \delta_{n,m},\\
& &[ P_i,~Q_j ] = A_{ij},
\end{eqnarray*}

\noindent where $(A_{ij})$ is the Cartan matrix for the Lie algebra $g$.
Let

\begin{eqnarray*}
& &A=B=1,\\
& &u[n]=\frac{(pq)^{-1/2}}{q^n-1},~~~~v[n]=\frac{q^{-1/2}}{(pq)^{-n}-1}
\end{eqnarray*}

\noindent then we have

\begin{prop}
The following bosonic expressions give a level $c=1$ realization for the algebra
\Egc{(1)}{q,p,c} with $\Psi_{ij}(z|q)$ chosen as in (\ref{dd}), 
on the Fock space of the Heisenberg algebra ${\cal H}_{q,p}(g)$,

\begin{eqnarray*}
& & E_i(z)= : \exp[ \varphi_i(z)] :,\\
& & F_i(z)= : \exp[ \psi_i(z)] :,\\
& & H^+_i(z)= z^{-2}:E_i(zp^{1/4}) F_i(zp^{-1/4}):,\\
& & H^-_i(z) =z^{-2}:E_i(zp^{-1/4}) F_i(zp^{1/4}):.
\end{eqnarray*}

\end{prop}

Of cause, this result can be readily mapped into a representation of 
\Egc{(2)}{q,\tilde{q},\beta,\gamma} at $\tilde{q}=qp, \beta=q, \gamma=qp$ using the 
homomorphism $\mu$.

\section*{}
\centerline{---------------------}
\vspace{0.2cm}

So far we obtained two families of quantum current algebras \{\Egc{(1)}{q_n,p,c_n}, $n\in Z$\}
and \{\Egc{(2)}{q^{(n)},\tilde{q}^{(n)},\beta_n,\gamma_n}, $n \in Z$\} and established
their structures as infinite Hopf family of algebras.
The generality of the defining relations for these two family of algebras 
indicates that the infinite Hopf family of algebras exists much broader than the standard Hopf
algebras. Actually, taken from the point of view of defining tensor product representations,
the standard Hopf algebra structure is by no means superior to the infinite Hopf family of algebras,
because both kinds of structures allow one to obtain fused representations from the tensor category
of the set of seed algebras. 

It is interesting to mention that the comultiplications appearing in
such co-structures are all of the Drinfeld type, which closes over the currents
themselves and does not require the resolution to the inverse problem (Riemann problem)
of the Ding-Frenkel homomorphism \cite{KLP1}.  
Recall that two kinds of comultiplications (and thus
two kinds of Hopf algebra structures) are known for the standard $q$ affine algebras.
While the $q$ affine algebras are considered as 
the most trivial cases of infinite Hopf family of algebras, 
only the Drinfeld type co-structures find their place in the generalized co-structure, 
whilst the standard Hopf algebra structure find no counterpart
in our present study. This is because we studied here only the current algebra 
formulation. In order to have a complete generalization of both Hopf algebra structures
of the $q$ affine algebras, it seems that we have to go to the Yang-Baxter realization 
as well, and it is highly probable that in that realization, the quasi-triangular
quasi-Hopf algebra structure may take some place. We leave this problem to future
study. 
 
We should emphasis that this work is only a preliminary study for the
new quantum current algebras. Besides the definition and infinite Hopf family structure,
we know very little about these algebras, especially their detailed representation
theory, vertex operators, Yang-Baxter type realizations etc. The physical applications
should also be considered.

On the other hand, the structures of infinite Hopf family of algebras is still poorly
understood yet.  We do not know whether there exists a quantum double construction over
the infinite Hopf family of algebras and, if not, what kind of new structure will take
the place of the standard quantum doubles. Also, the classical counterpart of the
infinite Hopf family of algebras is unknown and it seems that all these problems
deserve further investigations.

Finally, from experiences of studying various (deformed) affine algebras, we know that
given a (quantum deformed) affine algebra there must exist an accompanied (deformed)
Virasoro algebra, and the latter is highly expected to have important applications
in physics of 1+1 dimensions. Therefore, given the two new family of quantum current algebras, 
it seems very interesting to find/construct the corresponding deformed Virasoro algebras. 

\vspace{0.5cm}

{\bf Acknowledgement.} I would like to acknowledge the Abdus Salam ICTP for hospitality 
during the completion of this work, which was done under the Associateship Scheme of the
Abdus Salam ICTP, Trieste, Italy.

\end{document}